\documentclass[12pt]{article}

\usepackage{epsfig,amsfonts}
\usepackage{amsmath}
\usepackage{amssymb}

\usepackage{graphicx}

\renewcommand{\i}{\mathrm i}

\newcommand{\surf}[1]{\mathcal{#1}}
\newcommand{\grp}[1]{\mathcal{#1}}
\newcommand{\diag}[1]{\mathbf diag}
\renewcommand{\[}{\left[}
\renewcommand{\]}{\right]}
\newcommand{\complexC}{\mathbb{C}}

\renewcommand{\Re}{{\rm Re \,}}
\renewcommand{\Im}{{\rm Im \,}}

\newcommand{\be}{\begin{equation}}
\newcommand{\ee}{\end{equation}}

\newtheorem{theorem}{Theorem}
\newtheorem{prop}{Proposition}
\newtheorem{problem}{Problem}
\newtheorem{definition}{Definition}

\author{Andrey V. Shanin$^1$, Eugeny M. Doubravsky$^2$ \\
Russia, 119899, Moscow, Moscow State University, \\
$^1$ Department of Physics, $^2$Department of Mechanics and
Mathematics}

\title{Criteria for commutative factorization of a class of algebraic matrices}

\begin{document}

\maketitle

\begin{abstract}
The problem of matrix factorization motivated by diffraction or
elasticity is studied. A powerful tool for analyzing its solutions
is introduced, namely analytical continuation formulae are
derived. Necessary condition for commutative factorization is
found for a class of ``balanced'' matrices.
Together with Moiseyev's method and Hurd's idea, this gives a description
of the class of commutatively solvable matrices.
As a result, a simple analytical procedure is described, providing an answer,
whether a given matrix is commutatively factorizable or not.
\end{abstract}


\section{Introduction}

A matrix factorization problem (i.e.\ a problem of finding the factors
$Q^+$ and $Q^-$ providing the decomposition (\ref{problem_1}) for a known
matrix $G(k)$)
is usually motivated by an elasticity
or a wave diffraction problem. Typically its formulation does not contain
a requirement of {\em commutative\/} factorization. However, the
possibility to perform a commutative factorization is usually
studied carefully, since all known factorization methods are based on
ideas connected with commutativity.

There are two main methods for commutative factorization.
The first one is based on the idea
by Heins \cite{Heins48} who proposed to split the logarithm of matrix
$G$ additively. Further progress in this direction is connected with the
names of Chebotarev \cite{Chebotarev56} and Khrapkov \cite{Khrapkov71}.
In the paper by Khrapkov an explicit form of factorization for a certain
class of matrices $2 \times 2$ was found.

However, sometimes Khrapkov's formula leads to exponential
growth of the factors at infinity, and this is not acceptable for physical
applications. To suppress this growth in some cases,
a special technique was proposed
by Daniele \cite{Daniele84} and Hurd \& L\"uneburg \cite{Hurd85},
who thus enlarged the class of explicitly
factorizable matrices.
A factorization to some other matrices can be obtained if the analytical
continuation of the matrix is studied (see Rawlins \cite{Rawlins75} and
Hurd \cite{Hurd76}).
Matrices with dimension more than $2 \times 2$ have been investigated by
Lukyanov \cite{Lukyanov80}.

Another technique of commutative matrix factorization is based
on diagonalization of $G$ and studying the eigenvalues as a multi-valued
function.
This approach has been developed by Cercignani \cite{Cercignani77},
L\"uneburg \cite{Luneburg82}, Moiseyev \cite{Moiseyev89},
Meister \& Penzel \cite{Meister92} and
Antipov \& Silvestrov \cite{Antipov02}. It is based mainly
on the research of Zverovich \cite{Zverovich71}, who described a method for
solving scalar boundary-value problem on Riemann surfaces.

It was noticed by Hurd \cite{Hurd87} and Daniele
\cite{Daniele84a} that all known methods to factorize a matrix
are connected with the fact that
the matrix commutes with a polynomial matrix.

A question of finding a rational factor enabling to perform
a commutative factorization for matrices $2 \times 2$
was discussed by Williams \cite{Williams84}.
Also this question was stied in detais by Ehrhardt and Speck \cite{Speck02}.

Beside the exact methods, an interesting attempt to solve a matrix factorization
problem approximately was proposed by Abrahams \cite{Abrahams98}.

Each of the works having either Kharapkov's or Moiseyev's context was dedicated to a distinct
class of matrices, i.e.\ the starting point of such a work was a phrase like ``consider a matrix having
the following form \dots''. Sometimes, however, it is not easy to say whether a given matrix
can be reduced to one of the known classes by algebraic manipulations.
That is why,  a separate and interesting issue is a classification
of matrices with respect to factorization, i.e.\ finding a criteria, for example, of
the possibility of commutative factorization. Such a theory is known
for matrices $2 \times 2$ \cite{Chebotarev56}, however the
necessary condition for commutative
factorization has been found in an ``Ansatz'' form, i.e.\ a matrix
should have a specific representation including some entire
(polynomial) matrices and some arbitrary functions as
coefficients.

A necessary condition for commutative factorization for matrices of arbitrary dimension
was studied by Jones in \cite{Jones84}, however the author restricted himself to the
case of matrices having distinct eigenvectors {\em everywhere}, while a typical matrix
emerging in diffraction theory has branch points, i.e.\ it has distinct eigenvectors
{\em almost everywhere}.

The current paper is inspired mainly by the works of Antipov et.al.\ and the bright work of Hurd.
The idea is to take an algebraic matrix $G(k)$ and study the properties of
the factors $Q^+$ and $Q^-$ on their Riemann surfaces following {\em a priori}
from the decomposition (\ref{problem_1}), i.e.\ without constructing the factors
explicitly. We found that the decomposition (\ref{problem_1}) taken together with
the regularity conditions imposed on the factors define a unique
Riemann surface of $Q^+$ or $Q^-$. Moreover, the values of, for instance, $Q^+$
taken on different sheets are connected by simple algebraic relations.
Thus, in Section~2 we obtain {\em analytical continuation formulae}. The relation obtained
by Hurd is a particular case of such formulae.

The problem for $Q^+$ reminds the functional problem for Abelian integrals, but
while the value of an Abelian integral is increased by a constant when the argument
is carried along a certain contour, the value of $Q^+$ is multiplied by a
non-constant {\em bypass matrix}.

In Section~3 we study the question
of the possibility of commutative matrix factorization.
To describe the matrices, for which the necessary
condition for commutative factorization is fulfilled,
we introduce the class of {\em branch-commutative}
matrices $G(k)$, i.e.\ the matrices whose values on different sheets of their Riemann
surface commute. Branch-commutativeness is a much weaker condition than commutativeness
introduced by Chebotarev \cite{Chebotarev56}, since we demand the
commutation only of matrices having the
same affix~$k$. We demonstrate that branch-commutative matrices can be factorized
by Moiseyev's  technique \cite{Moiseyev89}, so the branch-commutativeness is necessary
and sufficient condition for
commutative factorization.

In Section~4 we study the Hurd's idea in our terms. {\em Bypass matrices} are introduced. They are
the matrices connecting the values of unknown function $Q^+$ on different sheets.
If all bypass matrices commute, we call matrix $G$ {\em bypass-commutative}. The class of
bypass-commutative matrices is wider than the class of branch-commutative, for example all algebraic
matrices, to which Hurd's method can be applied are bypass-commutative. We show that
for any bypass-commutative matrix $G$ a rational factor can be found, transforming $G$ into
a branch commutative matrix.

\section{Analytical continuation of the factors $Q^\pm$}

\subsection{Problem under consideration}

The initial problem of matrix factorization is as follows:

\begin{problem}
For a matrix $G(k)$ defined in a narrow strip along the real axis
($-\epsilon <\Im [k]< \epsilon$) find matrices $Q^+ (k)$, $Q^-
(k)$ analytical (maybe except some isolated poles), continuous,
having algebraic growth in the upper ($\Im [k] > -\epsilon$)
and lower ($\Im [k] < \epsilon$) half-planes, respectively, and
satisfying the equation
\begin{equation}
    G(k)=Q^+(k)Q^-(k).
\label{1.1.1}
\end{equation}
\label{problem_1}
\end{problem}

Algebraic growth hereafter means that there exists a number $l$,
such that all elements of corresponding matrices grow at the
corresponding half-plane no faster than $|k|^l$. We cannot expect
that the elements will grow exactly as some powers of $k$, since
the solutions can have logarithmic behaviour.

We assume that $G(k)$ is an algebraic function, as it happens in 
many applications. Thus, the function itself and the relation 
(\ref{1.1.1}) can be continued from the strip $-\epsilon <\Im [k]< \epsilon$. 

Besides, we assume everywhere that the determinant of $G$ is not equal to
zero identically.



\subsection{Notations for bypasses}

Let $\surf{R}_G$ be the Riemann surface of matrix $G(k)$. Below we shall call
$k$ an {\em affix\/} of a point $(k, G(k)) \in \surf{R}_G$.

Let branch points of $G(k)$ have affixes $\tau^+_m$ and $\tau^-_m$,
where $\Im[\tau^+_m] >0$ and $\Im[\tau^-_m]< 0$.
Make $G$ single-valued on $\mathbb{C}$ by performing cuts
going from branch points to infinity. The
cuts can be chosen as $\gamma^+_m =
(\tau_m^+, + \mathrm{i} \infty)$ and $\gamma^-_m = (\tau_m^-, - \mathrm{i} \infty)$.
It is
important that the cuts should not cross the real axis and each
other. As a result, the surface $\surf{R}_G$ becomes split into several
sheets. There is a special sheet of $\surf{R}_G$, on which equation
(\ref{1.1.1}) is assumed to be valid.  Name this sheet
a {\em physical sheet}.

Here and below, the structure of Riemann surface is displayed by
graphical diagrams. Horizontal lines correspond to the sheets,
nodes correspond to branch points,
and vertical lines link sheets, which are connected at a certain branch
point.

\begin{figure}[ht]
\centerline{\epsfig{file=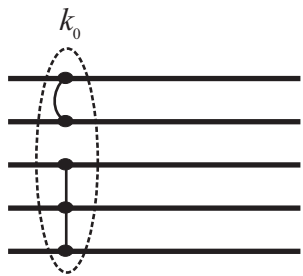}}
\caption{Order of an affix} \label{affix}
\end{figure}

In some artificial cases an affix may correspond to several branch points
having different orders. Define for each affix its order $n^{\pm}_m$, which is
the least common multiple of all orders of branch points with corresponding affix.
For example, in Fig.~\ref{affix} a fragment of a Riemann surface is shown.
Parameter $n$ for the affix $k_0$ is equal to~6.

Introduce a notation for the sheets of $\surf{R}_G$. Note that later the
same notation will be used for the sheets of the Riemann surfaces of $Q^{\pm}$.
Each point of the surface will be denoted by $(k)\{w\}$,
where $k$ is an affix, and $\{ w \}$ is a {\em word\/} describing
the path, along which the argument $k$ should be carried from
physical sheet to a selected sheet. The structure of this word is
explained below.

\begin{figure}[ht]
\centerline{\epsfig{file=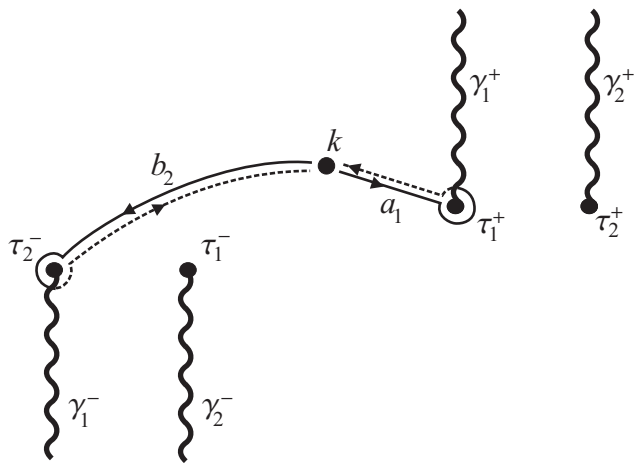}}
\caption{Notation for bypasses} \label{fig1}
\end{figure}

Denote bypasses about points
$\tau^+_i$ in positive direction by letters $a_i$ and bypasses about
points $\tau_i^-$ in positive direction by $b_i$ (Fig.~\ref{fig1}).

A series of consecutive bypasses will be denoted by a word of
letters $a_i$ and $b_i$. The word must be read from left to right,
i.e.\ the first performed bypass corresponds to the left end of
the word and the last bypass corresponds to the right end. By
default a series of bypasses begins from the ``physical sheet''. A
trivial bypass will be denoted by letter $e$.

Define the composition of words $w$ and $v$ as the bypass
performed along the way composed of $w$ and $v$. The bypass $w$ is
performed first. Denote this composition by $wv$. Let $\grp{W}$ be
the set of all words, and let $\grp{W}_a$, $\grp{W}_b$ be the sets
of words composed only of the letters $a_i$, and only of letters
$b_i$, respectively.

Let $G(k)\{e\}$ be the value of function $G(k)$ on the ``physical
sheet''. Denote by $G(k)\{w\}$ the value of $G(k)$ on the sheet
that can be reached by performing the bypass $w$ starting from the point
$(k,G(k)\{e\})$.

The set $\grp{W}$ can be considered as a {\em group\/} of words, a
subject of combinatorial group theory. Its generators are the
letters $a_i, b_i$, and the relations have the form
\be
a_i ^{n^+_i} = e, \qquad b_i ^{n^-_i} = e.
\label{branchpoints}
\ee
As we shall see below, the same relations are valid for the words describing
the Riemann surfaces of $Q^{\pm}$.

Relations (\ref{branchpoints}) enable one to determine an inverse element for each
$w \in \grp{W}$ without introducing new letters for bypasses in negative direction
(or without using the symbols $a_i^{-1}$ and $b_i^{-1}$).
Using (\ref{branchpoints}), below we assume that for each
word $w$ there exists a word $w^{-1}$, such that $w w^{-1} = w^{-1} w = e$.

Let us demonstrate an example of Riemann surfaces for $G$ and
$Q^+$. Take matrix $G$ from Daniele's paper \cite{Daniele84}:
\be
G(k) = \left( \begin{array}{cc}
1 & \frac{k_1 - s(k)}{k_2 + s(k)} \\
\frac{k_2 - s(k)}{k_1 + s(k)} & 1
\end{array} \right),
\label{test}
\ee
where $s(k) = \sqrt{k_0^2 - k^2}$; $k_0$, $k_1$
and $k_2$ are some complex constants.

In this case the Riemann surface of $G(k)$ has two sheets and two
quadratic branch points, namely $k = \pm k_0$. Let be $\Re[k_0]
>0$. Let letter $a$ denote a bypass about $k_0$, and letter $b$
denote a bypass about $-k_0$.

\begin{figure}[ht]
\centerline{\epsfig{file=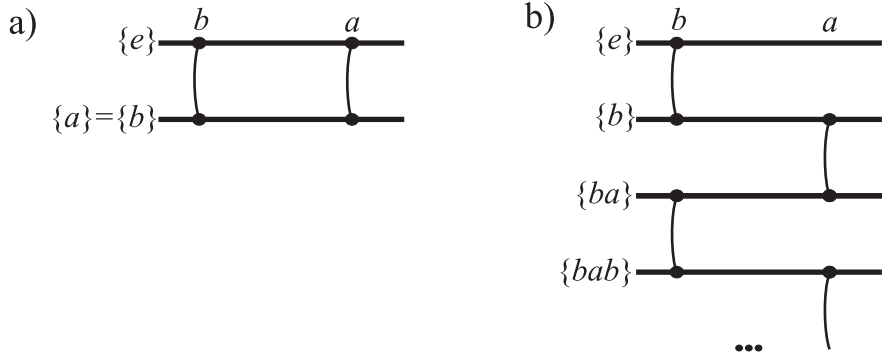}}
\caption{Diagrams of
Riemann surfaces for $G(k)$ and $Q^+ (k)$}
\label{diagrams}
\end{figure}

The scheme for the Riemann surface of $G$ is shown in
Fig.~\ref{diagrams}~a. The upper sheet is physical (i.e.\ it
contains the ``physical'' real axis).

The scheme of $Q^+$ corresponding to this problem is shown in
Fig.~\ref{diagrams}~b. The number of sheets is infinite, but all
branch points are of second order, and the positive physical half-plane
contains no branch points. This structure can be revealed, e.g.\ from
\cite{Daniele84}.

\subsection{Truncation operators}

Let be $w=\alpha_1\alpha_2\ldots\alpha_n$ where $\alpha_i$
substitutes an arbitrary single letter. Denote by $p$ the maximal
number, such that the word
$\alpha_1\alpha_2\ldots\alpha_p\in\grp{W}_a$. Analogically let $m$
be the maximal number, such that
$\alpha_1\alpha_2\ldots\alpha_m\in\grp{W}_b$. Obviously, one of
this integers is zero, since the first letter of the word is
either $a_j$ or $b_j$.

Define the truncation operators $^+$ and $^-$ by
\begin{eqnarray*}
&&w^+=\alpha_{p+1}\alpha_{p+2}\ldots\alpha_n,
\label{4.6}\\
&&w^-=\alpha_{m+1}\alpha_{p+2}\ldots\alpha_n.\ \label{4.7}
\end{eqnarray*}

For example, applying operators $^+$ and $^-$ to the words
$w=a_1a_2b_1b_2$, $v=b_1b_2a_1a_2$ we obtain
\begin{eqnarray*}
    & & w^+=b_1b_2, \qquad w^-=w=a_1a_2b_1b_2, \qquad w^{+-} \equiv (w^+)^- = e,
    \label{4.8}\\
    & & v^-=a_1a_2, \qquad v^+=v=b_1b_2a_1a_2, \qquad v^{-+}  = e.
    \label{4.9}
\end{eqnarray*}


\subsection{Formulae of analytical continuation}

Consider equation (\ref{1.1.1}). Both right
and left sides of this equation are analytic functions in some
neighbourhood of the real axis of the physical sheet. Continue
$Q^+$ and $Q^-$ analytically to this domain and, further, onto
some Riemann surfaces. Continue also the relation (\ref{1.1.1})
onto the Riemann surfaces of $G$, $Q^+$ and $Q^-$. Obviously, the
continuation of the relation (\ref{1.1.1}) can be written in the
form:
\be
Q^+(k)\{w\} \, Q^-(k)\{w\}=G(k)\{w\}.
\label{2.1.1}
\ee
At this formula (\ref{2.1.1})
has sense only for geometrically fixed bypasses.

Here we are going to find the {\em formulae of analytical
continuation\/} for $Q^{\pm}$, i.e.\ algebraic relations
connecting $Q^{\pm} (k)  \{ w \}$ with $ Q^{\pm} (k) \{ e \} $.

General formulae of analytical continuation can
be written in a recursive form as follows:

\begin{theorem}
Let $Q^+ (k)$ and $Q^- (k)$ form a solution of
Problem~\ref{problem_1}. Then the following relations are valid
\be
Q^+\{w\}=G\{w^+\} \, G^{-1} \{w^{+-}\} \, Q^+\{w^{+-}\} ,
\label{2.1.3}
\ee
\be
Q^-\{w\}= Q^-\{w^{-+}\} \, G^{-1} \{w^{-+}\} \, G\{w^-\} .
\label{2.1.4}
\ee
(A dependence on $k$ is
implied for all functions in (\ref{2.1.3}), (\ref{2.1.4})).
\label{prop_2.1}
\end{theorem}

The proof is rather straightforward and based on the relations
following from the regularity conditions
\begin{equation}
Q^+\{w\}=Q^+\{w^+\}, \qquad Q^-\{w\}=Q^-\{w^-\}
\label{4.10}
\end{equation}
literally denoting that $Q^+$ is analytical at the points $\tau^+_j$, while 
$Q^-$ is analytical at the points $\tau^-_j$.
Let us prove (\ref{2.1.3}) (relation (\ref{2.1.4}) is similar).
First, according to (\ref{4.10}),
\begin{equation}
Q^+ \{  w \} = Q^+ \{ w^+ \}.
\end{equation}
Then, according to (\ref{2.1.1}),
\begin{equation}
Q^+ \{  w \} = Q^+ \{ w^+ \} = G \{ w^+ \} (Q^- \{ w^+ \} )^{-1}. 
\end{equation}
According to the second relation of (\ref{4.10}),
\begin{equation}
Q^+ \{  w \} = Q^+ \{ w^+ \} = G \{ w^+ \} (Q^- \{ w^+ \} )^{-1} = G \{ w^+ \} (Q^- \{ w^{+-} \} )^{-1}.
\end{equation}
Finally, according to (\ref{2.1.1})
\begin{equation}
(Q^- \{ w^{+-} \} )^{-1} = (G\{w^{+-}\})^{-1} Q^+ \{ w^{+-} \}
\end{equation}
and we get (\ref{2.1.3}).

Note that for any word $w$ their exists some constant $c$, such
that $w^{(+-)^c} = e$, therefore  formula (\ref{2.1.3}) being repeated
several times
connects $Q^+\{ w \}$ with $Q^+ \{ e\}$. Analogously,
$Q^-\{ w \}$ is connected with $Q^- \{ e \}$. The coefficients are
always products of known matrices.

Analytical continuation in the form (\ref{2.1.3}) has been
obtained by Hurd \cite {Hurd76} for a particular case of a single
bypass. Hurd's ideas are discussed later in details.

Using analytical continuation we can investigate the structure of
Riemann surface of unknown function $Q^+$. For example, the
following proposition can be easily proved:

\begin{prop}
Let $G(k)$ be an algebraic matrix, and let the functions $Q^+
(k)$ and $Q^- (k)$ form  a solution of Problem~1. Then the
functions $Q^+ (k)$ and $Q^- (k)$ can be analytically continued
onto some Riemann surfaces; both functions have branch points only
at affixes $\tau_i^{\pm}$. The order of each branch point is a divisor
of corresponding $n^{\pm}_i$.
\end{prop}

A formal proof can be conducted by induction with respect to the
length of the word $w$, which is the argument of $Q^{\pm}(k)\{ w \}$.

Generally, solution of Problem~\ref{problem_1} is not unique: for
example the behaviour of different solutions at infinity can be
different. However, it is easy to prove that all solutions are similar
up to a meromorphic matrix factor.

\section{Necessary condition for commutative matrix factorization}

\subsection{Necessary condition in the ``check-up'' form}

\begin{definition}
Let $G(k)$ be an algebraic matrix, let its branch points have affixes $\tau^{\pm}_j$
and lie aside from the real axis. Let the sets $\grp{W}$, $\grp{W}_a$, and $\grp{W}_b$
be defined as described above.
Riemann surface $\surf{R}_G$ will be called balanced if for any $w \in \grp{W}$ there exist
words $w_a \in \grp{W}_a$ and  $w_b \in \grp{W}_b$ such that
\be
G(k)\{w_a\} = G(k)\{w_b\} = G(k)\{w\}.
\label{def_bal}
\ee
\end{definition}

An example of a balanced Riemann surface is the Riemann surface of scalar function
$\sqrt{1+\sqrt{2+k^2}}$ with an arbitrary choice of the physical sheet. Besides, a surface of
an arbitrary matrix function, which is a rational combination of $k$ and several square roots
$\sqrt{\tau^2_j-k^2}$, is balanced.

An example of a Riemann surface that is not balanced is the surface of the
function $\sqrt{\i + k} + \sqrt{-\i+k}$.

\begin{definition}
Algebraical matrix $G(k)$ is called branch-commutative, if for any $k$ the
values of $G$ on different sheets of its Riemann surface commute, i.e.
\be
[ G(k)\{ w_1 \} , G(k)\{ w_2 \} ] \equiv
G(k)\{ w_1 \} G(k)\{ w_2 \} - G(k)\{ w_2 \}  G(k)\{ w_1 \} = 0.
\label{def_com}
\ee
for any different words $w_1$ and $w_2$.
\end{definition}

To illustrate the definition of branch-commutativeness
consider a simple example.
Take a matrix
\be
G(k) = \left( \begin{array}{ccc}
k & 2k &  s(k) \\
2k & k & -s(k) \\
-s(k) & s(k) & k
\end{array} \right), \qquad s(k) = \sqrt{k_0^2-k^2}.
\label{example}
\ee
Let be $\Im[k_0] >0$. There are two letters, $a$ and $b$, corresponding to bypasses
about $k_0$ and $-k_0$.  The value of $G$
on the sheet $\{a\}$ is equal to
$$
G(k)\{a\} = \left( \begin{array}{ccc}
k & 2k &  -s(k) \\
2k & k & s(k) \\
s(k) & -s(k) & k
\end{array} \right).
$$
Since $\surf{R}_G$ has two sheets, to check branch-commutativeness
one should check only the identity
\be
[G(k)\{e\}, G(k)\{a\}] =0,
\label{example1}
\ee
where $G(k)\{e\}$ is defined by (\ref{example}).
Simple calculations show that (\ref{example1}) is fulfilled, so (\ref{example})
is a branch-commutative matrix.

The necessary condition of commutative factorization is given by
the following theorem:

\begin{theorem}
If a matrix $G$ having balanced Riemann surface admits commutative factorization
\be
Q^+(k) \, Q^-(k) = Q^-(k) \, Q^+(k) = G(k),
\label{3.2.1}
\ee
then it is a branch-commutative matrix.
\label{prop_main}
\end{theorem}

\noindent {\bf Proof: } The formula of analytical continuation
(\ref{2.1.3}) has been derived for right
factorization. One can obtain a similar formula for left
factorization $G(k)=Q^-(k) \, Q^+(k)$:
\be
Q^+\{w\}=Q^+\{w^{+-}\} \, G^{-1}\{w^{+-}\}
G\{w^{+}\}.
\label{3.2.2}
\ee

Perform the rest of the proof step by step. Here we mark the
statements and make some comments if the statements are not
obvious:
\begin{enumerate}

\item
For any word $w$
$ \quad
Q^+\{w\}Q^-\{w\}=Q^-\{w\}Q^+\{w\}=G\{w\}        
$.
It is an analytical continuation of (\ref{3.2.1}).

\item
For any word $w$
$ \quad
\[G\{w\},(Q^+\{w\})^{-1}\]=0
$.
This follows from $G\{ w \} (Q^+\{ w \})^{-1} = Q^-\{ w \} = (Q^+\{ w \})^{-1} G\{ w \}$. 

\item
For any word $w$
$ \quad
\[G\{w\},Q^+\{w\}\]=0
$.
This can be obtained from the
previous point by multiplication by $G^{-1}$ at left and right.

\item
For any $v\in\grp{W}_a$
$ \quad
\[G\{v\},Q^+\{e\}\]=0
$. This follows from the previous point and (\ref{4.10}).

\item
For any word $w$
$ \quad
\[G\{w\},Q^+\{e\}\]=0
$.
Note that for a matrix with balanced Riemann surface
for any word $w$ there exists a word $v \in \grp{W}_a$, such that
$G\{ v\} = G\{w\}$, and $Q^+\{v\} = Q^+\{e\}$.

\item
For any $v_1 \in \grp{W}_b$, $v_2 \in \grp{W}_a$
$ \quad
\[ G \{ v_1 v_2 \} , G^{-1} \{  v_2 \} \] =0
$. This statement
can be obtained by applying left and right analytical continuation
formulae to the word $v_1 v_2$ and by using the previous point.

\item
For any $v_1 \in \grp{W}_b$, $v_2 \in \grp{W}_a$
$ \quad
\[ G \{ v_1 v_2 \} , G \{  v_2 \} \] =0
$.

\item
The statement of the theorem, by noting that for any $w_1$ and
$w_2$ one can find the words
$v_b \in \grp{W}_b$ and $v_a \in \grp{W}_a$,
such that $G\{ v_b v_a \} = G\{ w_1\} $,  $G\{ v_a \} = G\{ w_2\} $.

\end{enumerate}

Theorem~\ref{prop_main} is an important result of the paper. Note that
since the number of sheets of $G$ is finite, the necessary
condition can be established by checking a finite number of
matrix identities.


\subsection{Diagonalization and properties of eigenvectors}

Let an algebraic (not necessarily branch-commutative)
matrix $G(k)$ have distinct eigenvalues
almost everywhere (i.e.\ on the whole complex plane excluding
several points). Represent this matrix in the form
\be
G(k) = M(k) \mathop{\rm diag} \[ \lambda_1, \dots ,\lambda_N \] M(k)^{-1}
\label{4.1.1}
\ee
Here matrix $M(k)$ consists of vector-columns,
which are right eigenvectors of $G$; $\lambda_1 \dots \lambda_N$
are  corresponding eigenvalues; $N$ is dimension of $G$.
Normalize the columns of $M$ by making all elements of the first
raw of $M$ equal to~1.

Obviously, for obtaining representation (\ref{4.1.1}) one
should first solve the characteristic equation for $G$, and then find
a solution of an inhomogeneous linear system for each eigenvector.

Denote Riemann surface of matrix $M(k)$ by $\surf{R}_M$.
Now we have associated with a matrix $G$ two Riemann surfaces:
$\surf{R}_G$ and $\surf{R}_M$. Typically, say for Khrapkov
matrices, $\surf{R}_M$ has a structure very different from
$\surf{R}_G$.

A lot of authors studied matrix
factorization problems by formulating a functional problem on a
Riemann surface. It is important to mention that most of them
had in mind the surface $\surf{R}_M$, not $\surf{R}_G$.
Typically, the surface $\surf{R}_M$ is studied in Moiseyev's context,
and $\surf{R}_G$ in Hurd's one.

Obviously, Riemann surface for the eigenvalues $\lambda_j (k)$
should contain branch points of both structures, i.e.\ of $G$ and
of $M$.

Let $G(k)$ be a {\em branch-commutative\/} matrix. In this case the set
of normalized eigenvectors must be the same on all sheets of
$\surf{R}_M$. Therefore, matrix $M(k)$ possesses an important
property: any bypass about branch points leads to a permutation of
the columns, i.e.\ an analytical continuation of each column along
a closed contour $c$ on $\mathbb{C}$ is some other column of~$M$.

\begin{figure}[ht]
\centerline{\epsfig{file=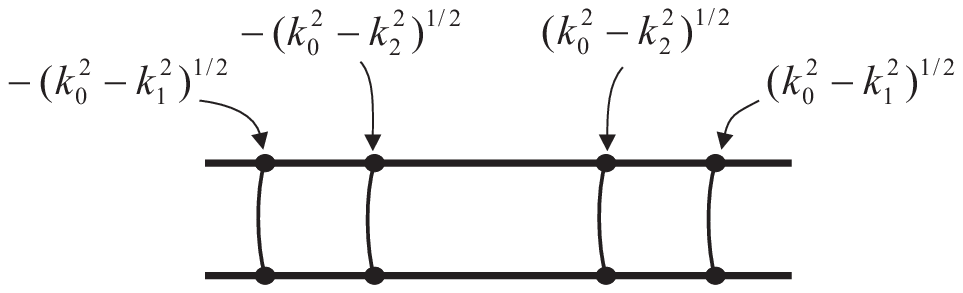}}
\caption{Diagram of $\surf{R}_M$}
\label{diagram_M}
\end{figure}

As an example, consider matrix (\ref{test}), which is
branch-commutative. As it was mentioned, it has only two branch
points, namely $\pm k_0$. The scheme of Riemann surface for this
matrix is shown is Fig.~\ref{diagrams}~a. It is easy to find that
matrix $M$ for this $G$ is as follows:
\be
M(k) = \left(
\begin{array}{cc}
1 & 1 \\
\frac{\sqrt{k_0^2 - k^2 - k_2^2}}{\sqrt{k_0^2 - k^2 - k_1^2}} & -
\frac{\sqrt{k_0^2 - k^2 - k_2^2}}{\sqrt{k_0^2 - k^2 - k_1^2}}
\end{array} \right).
\label{4.1.3}
\ee
Matrix $M$ has four branch points, namely
$\pm \sqrt{k_0^2-k_1^2}$ and $\pm \sqrt{k_0^2-k_2^2}$. Generally (i.e.\
if $k_1 \ne 0$ and $k_2 \ne 0$) the branch points of $\surf{R}_M$
are different from the branch points of $\surf{R}_G$. The scheme
of $\surf{R}_M$  is shown in Fig.~\ref{diagram_M}.

A transition from one sheet of $\surf{R}_M$ to another leads
to a permutation of the columns of~$M$.


\subsection{``Ansatz'' form of necessary condition}

\label{Ansatz}

\begin{theorem}
Let $G$ be a branch-commutative matrix $N \times N$, whose eigenvalues
are distinct almost everywhere. Then it can be represented in the
form
\be
G = \sum_{m=0}^{N-1} g_m( k ) A^{m} (k),
\label{4.2.1}
\ee
where $A (k)$ is a rational matrix, $g_m (k)$ are algebraic
functions. Vice versa, any matrix admitting a decomposition of the
form (\ref{4.2.1}) is branch-commutative.
\label{prop_4.1}
\end{theorem}

\noindent
{\bf Proof: } The second part of the theorem is obvious, so we are
concentrating our efforts on the first one. Consider matrix
$M(k)$. Let $\pi_c$ be a permutation of columns of $M$ occuring
when the argument is carried along a contour $c$ on $\mathbb{C}$
starting and terminating at $k$. Let $\Pi_c$ be a matrix
containing only numbers 0 and 1, describing permutation $\pi_c$ in
matrix language, i.e.\
\be
(\Pi_c)^m_n = \delta_{m, \pi_c (n)},
\label{star}
\ee
and the permutation of columns of $M$ looks like $M \to M \Pi_c$.

Construct $N$ functions $f_m (k)$, $m = 1 \dots N$ as follows.
Take $N$ constants $\beta_1 \dots \beta_N$ such that the
combinations
\be
f_m(k) = \sum_{n=1}^N \beta_n (M)^n_m ,
\label{4.2.3}
\ee
 almost everywhere obey the relation $f_{m_1}(k) \ne f_{m_2}(k)$ as
$m_1 \ne m_2$. (Here $(M)^n_m $ are the elements of $M$.)
Obviously, $f_m \to f_{\pi_c (m)}$ when the argument is carried along $c$.

Construct a combination
\be
A(k) = M(k) \mathop{\rm diag} \[f_1 (k), \dots , f_N (k) \]  M^{-1}(k).
\label{4.2.2}
\ee
Note that the diagonal matrix obeys the relation
\be
\mathop{\rm diag} \[ f_{\pi_{c}(1)} , \dots , f_{\pi_c (N)} \] =
\Pi_c^{-1} \mathop{\rm diag} \[f_1 , \dots , f_N  \] \Pi_c.
\label{starstar}
\ee
Substituting (\ref{star}) and (\ref{starstar}) into (\ref{4.2.2}),
conclude that $A$ remains unchanged after any bypass~$c$.
Since $A$ is an algebraic matrix by construction,
it should be a rational matrix.

Finally,
let us show that $G$ can be expressed in the form (\ref{4.2.1})
with matrix $A$ constructed above. The matrix composed of the
elements $(F)^m_n = f_n^{m-1} $ (here $m-1$ is a power, $m=1 \dots N$)
has a
non-zero determinant almost everywhere. In the opposite case it
would happen that $N$ distinct numbers are roots of a polynomial
of order smaller than $N$. Therefore any set of $N$ numbers, for
example the eigenvalues of $G$, can be represented as
\be
\lambda_n (k)  = \sum_{m=1}^N g_m (k) f^{m-1}_n (k)
\label{10.1}
\ee
for almost all $k$. By construction, $g_n$ are algebraic
functions.

The theorem is proved.

The form (\ref{4.2.1}) is close to that of \cite{Jones84}, however
on one hand we impose no restrictions on the behaviour of the
matrices $Q^{\pm}$, and on the other hand, we do not specify the
form of equation, which matrix $A$ should obey.

Theorem~\ref{prop_4.1} states that there are two alternative ways
to check, whether a diffraction matrix $G$ can be factorized
commutatively: 1) by checking whether a matrix can be
represented in a certain form, and 2) by checking conditions
(\ref{def_com}) between different sheets. The second variant seems
more easy.

\noindent
{\bf Note:} Commutative factorization of the matrices having form
(\ref{4.2.1}) was considered in \cite{Moiseyev89}, where an explicit
formula for the factors was constructed.

Thus, branch-commutativeness
(or, alternatively, the form (\ref{4.2.1}))
is a necessary and sufficient condition form commutative factorization
of matrices with balanced Riemann surfaces.


\section{Bypass matrices and Hurd's method}


\subsection{Bypass matrices}

Let $G$ be an algebraic matrix with a balanced Riemann surface,
and let $\grp{W}$ be the set of words associated
with this matrix. Let the right factorization problem (\ref{1.1.1}) be studied.

\begin{definition}
A bypass matrix $P_w (k)$ for a word $w$ is defined by the relation:
\be
P_w (k)=Q^+ (k) \{ w \} \, (Q^+ (k) \{ e \})^{-1}  .
\label{def_P1}
\ee
\end{definition}

\noindent
According to the formulae of analytical continuation
(\ref{2.1.3}),
\be
P_w (k) = (G\{w^+\} G^{-1}\{w^{+-}\})
(G\{w^{+-+}\} G^{-1}\{w^{+-+-}\}) \ldots
\label{def_P}
\ee
The product in the r.-h.s.\ is finite, since after some truncations the matrices
become equal to $G\{e\}$.


Let the number of sheets of $\surf{R}_G$ be equal to $n$.
Among all bypass matrices we can select a finite set of $n-1$ basic bypass matrices
\be
\hat P_j (k)= G (k)\{ w_j \} G^{-1} (k)\{ e \}
\label{bpm_1}
\ee
where $w_1 \dots w_{n-1}$ are any words belonging to $\grp{W}_b$, such that all
$G (k)\{ w_j \}$ belong to different sheets of $\surf{R}_G$, and none of these
sheets is the physical one.
Any bypass matrix can be written as a product of several basic bypass matrices
$\hat P_j$ taken in positive or negative  powers.


\subsection{Hurd's idea and its formalization}

The idea of Hurd \cite{Hurd76} can be expressed as follows: sometimes
the bypass matrices $P_w (k)$ can have a structure simplier than that of
the matrix $G(k)$. Here we express this simplicity in the following form.

\begin{definition}
Let $G(k)$ be an algebraic matrix with a balanced Riemann surface.
Let $P_w$ be a set of corresponding bypass matrices.
Let all basic bypass matrices commute with each other:
\be
[\hat P_{j} (k),  \hat P_{m} (k) ] = 0, \qquad j,m=1 \dots n-1,
\label{6.1.1}
\ee
The matrix $G$ will be called bypass-commutative.
\end{definition}

The approach introduced by Hurd is closely connected with bypass-commutativity. His idea was to
study the matrix boundary value problem on cuts made from branch points to infinity. The
matrix coefficient for each cut is one of the bypass matrices. All known techniques available
for such problems at the current moment require commutation of these bypass matrices and their analytical
continuations. That is why all known matrices, to which Hurd's method was applied successfully
are bypass-commutative.

The class of bypass-commutative matrices is quite wide. For example, all matrices $G$ with hyperelliptic
Riemann surface are bypass-commutative, since there is a single basic bypass matrix $G\{a\} G^{-1}\{e\}$.
An example of bypass-commutative matrix with a more sophisticated Riemann surface is
\be
G(k) = \left(  \begin{array}{cc}
s_1(k)   &   s_2(k) \\
-s_2(k)   & k s_1(k)
\end{array}  \right),
\label{exmpl}
\ee
where $s_1 = \sqrt{k_1^2-k^2}$, $s_2 =\sqrt{k_2^2- \sqrt{k_1^2-k^2}}$, $k_1$ and $k_2$
are constants.
The Riemann surface for such a matrix has four sheets; affixes of branch points are
$\pm k_1$ and $\pm \sqrt{k_1^2-k_2^2}$. Thus, there are three basic bypass matrices,
and their commutativity can be checked explicitly.

\begin{theorem}
Let matrix $G(k)$ be bypass-commutative.
Then there exists a rational matrix $S(k)$ such that the matrix
$G(k) S(k)$ is branch-commutative.
\label{th_hurd}
\end{theorem}

\noindent
{\bf Proof. }
Obviously, if basic bypass matrices commute then all bypass matrices commute.

Let $\surf{R}_G$ has $n$ sheets, and the dimension of $G$ be $N \times N$.

The matrix $S(k)$ can be constructed as follows:
\be
S(k) = \sum_{j=1}^{n} f(k) \{w_j \} G^{-1} (k) \{ w_j\},
\label{20.1.1}
\ee
where $w_j$ is a set of words listing all sheets of $\surf{R}_G$. For example,
the words from (\ref{bpm_1}) can be taken as $w_j$ for $j=1\dots n-1$, and $w_n = e$.

Function $f(k)$ is an arbitrary function, such that it is single-valued on $\surf{R}_G$, and
the r.h-s.\ of (\ref{20.1.1}) has non-zero determinant almost everywhere.
As a possible choice, one can construct $f(k)$ by the formula
\be
f(k) = \beta_{0,0} + \sum_{i=1}^N \sum_{j=1}^N \beta_{i,j} (G^{-1}(k))^i_j
\label{20.1.2}
\ee
with almost arbitrary constants $\beta_{i,j}$.

Since the sum in (\ref{20.1.1}) is taken over all sheets of an algebraic function $fG$,
the result is single-valued on $\complexC$, and therefore it is a rational matrix function.

By construction, the combination $G(k) S(k)$ is a linear combination of some bypass matrices
or their products. The analytical continuations of $G(k) S(k)$ can also be represented as products
of bypass matrices. Since the bypass matrices commute with
each other, $G(k) S(k)$ is a branch-commutative matrix.

An example (maybe quite simple) of such a consideration can be constructed using matrix (\ref{exmpl}).
Take function $f(k) \equiv 1$. Matrix $S$ is constructed by summation of $G^{-1}$ over four sheets of
its Riemann surface:
\be
S(k) = \frac{4(k_1^2 -k^2)}{k_2^4+k_1^4 k^2+(k^2-k_1^2)(1-2k_2^2 k + 2k^4)}
\left( \begin{array}{cc}
k & 0 \\
0 & 1
\end{array} \right)
\label{exmpl1}
\ee
Simple calculations show that $G S$ is a branch-commutative matrix.

Note that Moiseyev's approach can be applied to $GS$, i.e.\ an explicit commutative factorization
can be constructed.

Note also that if for some matrix $G$ there exists a rational matrix $S$, such that $G S$ or $S G$
is branch-commutative, then obviously all bypass matrices of $G$ commute. That is why,
bypass-commutativity is a necessary and sufficient condition for a balanced algebraic matrix
of a possibility to be converted into a branch-commutative matrix.


\section{A short summary}

The main results of this paper are as follows:

\begin{enumerate}

\item
Formulae of analytical continuation  are derived.

\item
Necessary condition of commutative factorization is found.
Namely, a balanced algebraic
matrix should be branch-commutative. This property can be
easily checked.

\item
Connection with the ``Ansatz'' form of the necessary
condition is established.

\item
Hurd's method is formalized. It can be applied if a matrix is bypass-commutative.
It is shown that in this case one can reduce the problem to
the commutative factorization case by multiplication by a rational matrix.

The necessary condition of commutative matrix factorization for balanced matrix
is checked as follows.
First, one should check, whether the values $G(k)$ taken on different sheets commute.
If they commute, then the matrix can be factorized by Moiseyev's method. Second,
one should construct the basic bypass matrices $P_j$ and check whether they commute
with each other. If they commute, then there exists a rational matrix $S$, multiplication by
which transforms matrix $G$ into a commutatively factorizable case.

\end{enumerate}

\end{document}